\documentclass[11pt]{article}
\usepackage{amsfonts}
\usepackage{amsmath}
\usepackage{amssymb}
\begin{document}

\baselineskip 18pt
\def\o{\over}
\def\e{\varepsilon}
\title{\Large\bf  Shifted\ Character\ Sums\ with\ Multiplicative\
Coefficients, II}
\author{K.\ Gong,\ C.\ Jia\ and\ M. A.\ Korolev}
\date{}
\maketitle {\small \noindent {\bf Abstract.} Let $f(n)$ be a
multiplicative function with $|f(n)|\leq 1,\,q$ be a prime number
and $a$ be an integer with $(a,\,q)=1,\,\chi$ be a non-principal
Dirichlet character modulo $q$. Let $\varepsilon$ be a sufficiently
small positive constant, $A$ be a large constant, $q^{\frac12+
\varepsilon}\ll N\ll q^A$. In this paper, we shall prove that
$$
\sum_{n\leq N}f(n)\chi(n+a)\ll N\frac{\log\log q}{\log q}
$$
and that
$$
\sum_{n\leq N}f(n)\chi(n+a_1)\cdots\chi(n+a_t)\ll N\frac{\log\log
q}{\log q},
$$
where $t\geq 2,\,a_1,\,\ldots,\,a_t$ are distinct integers modulo
$q$.}

\vskip.3in
\noindent{\bf 1. Introduction}

Let $q$ be a prime number, $a$ be an integer with $(a,\,q)=1,\,\chi$
be a non-principal Dirichlet character modulo $q$.

Since the 1930s, I.M.~Vinogradov had begun the study on character
sums over shifted primes
$$
\sum_{p\leq N}\chi(p+a),
$$
and obtained some deep results \cite{Vinogradov_1938},
\cite{Vinogradov_1953} where $p$ runs through prime numbers. His
best known result is a nontrivial estimate for the range $q^{{3\o
4}+\e}\ll N\ll q^A$, where $\e$ is a sufficiently small positive
constant and $A$ is a large constant. Later, A.A.~Karatsuba
\cite{Karatsuba_1970} widen the range to $q^{{1\o 2}+\e}\ll N\ll
q^A$.

For the M\"obius function $\mu(n)$, one can get same results on
sums
$$
\sum_{n\leq N}\mu(n)\chi(n+a)
$$
as that on sums over shifted primes. K.~Gong and C.~Jia
\cite{Gong_Jia_2015} considered the general sum
$$
\sum_{n\leq N}f(n)\chi(n+a), \eqno (1.1)
$$
where $f(n)$ is a multiplicative function with $|f(n)|\leq 1$. They
applied a modification (see \cite{Gong_Jia_2016}) of the method in
Section 2 in \cite{Bourgain_Sarnak_Ziegler_2013}, which is called as
the finite version of Vinogradov's inequality, to give a nontrivial
estimate for the sum (1.1) when $q$ is in a suitable range.
Precisely, they proved that if $f(n)$ is a multiplicative function
with $|f(n)|\leq 1,\,q$ is a prime number and $a$ is an integer with
$(a,\,q)=1,\,\chi$ is a non-principal Dirichlet character modulo
$q,\,q^{1\o 2}\leq N$, then
$$
\sum_{n\leq N}f(n)\chi(n+a) \ll {N\o q^{1\o 4}}\log\log(6N)+q^{1\o
4}N^{1\o 2}\log(6N)+{N\o \sqrt{\log\log(6N)}},
$$
which is nontrivial for $q^{{1\o 2}+\e}\ll N\ll q^A$.

In this paper, we shall use some new ideas from \cite{Korolev_2016b}
to give a refinement on the above upper bound.

{\bf Theorem}. \emph{Let $f(n)$ be a multiplicative function with
$|f(n)|\leq 1,\,q$ be a prime number and $a$ be an integer with
$(a,\,q)=1,\,\chi$ be a non-principal Dirichlet character modulo
$q$. Let $\e$ be a sufficiently small positive constant, $A$ be a
large constant, $q^{{1\o 2}+\e}\ll N\ll q^A$. Then we have}
$$
\sum_{n\leq N}f(n)\chi(n+a)\ll N{\log\log q\o \log q}, \eqno (1.2)
$$
\emph{where the implied constant depends on} $\e,\,A$.

\emph{For $t\geq 2$ and distinct integers $a_1,\,\ldots,\,a_t$
modulo $q$, we have}
$$
\sum_{n\leq N}f(n)\chi(n+a_1)\cdots\chi(n+a_t)\ll N{\log\log q\o
\log q}, \eqno (1.3)
$$
\emph{where implied constant depends on} $\e,\,A$ \emph{and} $t$.

Throughout this paper, we assume that $q$ is a sufficiently large
prime number, $\e$ is a sufficiently small positive constant, $A$
is a large constant. Let $p$ denote prime number. Set
$$
X=\log^2 q,\qquad\quad Y=q^{\e\o 4}. \eqno (1.4)
$$

\noindent{\bf 2. The proof of Theorem}

We need one auxiliary assertion. Write
\begin{align*}
R&=\{n:\,n\leq N,\, n{\rm\ has\ a\ prime\ factor\ in\ }(X,\,
Y]\},\\
T&=\{n:\,n\leq N,\, n{\rm\ has\ no\ prime\ factor\ in\ }(X,\,Y]\}.
\end{align*}

{\bf Lemma 1}. \emph{We have}
$$
|T|\ll N{\log\log q\o \log q}.
$$

{\bf Proof.} Let $\Phi(x,y)$ be the number of $n\le x$ free of prime
divisors $\le y$ and denote by $\Psi(x,y)$ the number of $n\le x$
free of prime divisors $>y$. Then, uniformly in $2\le y\le x$, we
have
\[
\Phi(x,y)\,\ll\,\frac{x}{\log{y}},\quad \Psi(x,y)\,\ll\,xe^{-{u\o
2}},\quad u\,=\,\frac{\log{x}}{\log{y}}. \eqno(2.1)
\]

The proof of these inequalities can be found in \cite{Korolev_2016a}
and in \cite[Part III, Ch. III.5, \S 5.1]{Tenenbaum_1995}
correspondingly.

Any $n\in T$ has the form $uv$: (a) all prime divisors of $u$ do not
exceed $X$ (or $u=1$), and all prime divisors of $v$ are greater
than $Y$; (b) all prime divisors of $u$ do not exceed $X$ (or $u=1$)
and $v=1$. If $N_{a}$ and $N_{b}$ denote the number of $n\in T$
satisfying (a) and (b) respectively, then $|T|\le N_{a}+N_{b}$.

Fixing the factor $u$, we obtain at most $\Phi(Nu^{-1},Y)$
possibilities to choose the factor $v = nu^{-1}\le Nu^{-1}$. If
$v\ne 1$, then $v>Y$ and hence $Nu^{-1}>Y$. By (2.1), we get
$$
\Phi\biggl(\frac{N}{u},Y\biggr)\ll\frac{N}{\log{Y}}\cdot
\frac{1}{u}.
$$
The summation over $u$ yields
\begin{align*}
N_{a}&\ll\frac{N}{\log{Y}}\sum_{u}\frac{1}{u}\ll\frac{N}{\log{Y}}
\prod_{p\le X}\biggl(1+\frac{1}{p}+\frac{1}{p^{2}}+\cdots\biggr)\\
&\ll\frac{N}{\log{Y}}\prod_{p\le X}\biggl(1-\frac{1}{p}\biggr)^{-1}
\ll N\frac{\log{X}}{\log{Y}}\ll N \frac{\log\log{q}}{\log{q}}.
\end{align*}
Using (2.1) again, we obtain
$$
N_{b}\le \Psi(N,X)\ll N\exp{\biggl(-\frac{\log{N}}{2\log{X}}\biggr)}
\ll N\frac{\log{X}}{\log{N}}\ll N\frac{\log\log{q}}{\log{q}}.
$$
Hence, Lemma 1 is proved.

Now we begin the proof of Theorem. The discussion in Lemma 1 of
\cite{Gong_Jia_2015} yields
$$
\sum_{n\leq N}f(n)\chi(n+a)=\sum_{\substack{n\leq N\\
(n,\,q)=1}}f(n)\chi(n+a)+O\Bigl({N\o q}+1\Bigr).
$$
By Lemma 1,
$$
\sum_{\substack{n\leq N\\ (n,\,q)=1}}f(n)\chi(n+a)=\sum_{\substack{n\leq N\\
n\in R\\ (n,\,q)=1}}f(n)\chi(n+a) +O\Bigl(N{\log\log q\o \log
q}\Bigr).
$$
Therefore we have
$$
\sum_{n\leq N}f(n)\chi(n+a)=\sum_{\substack{n\leq N\\
n\in R\\ (n,\,q)=1}}f(n)\chi(n+a)+O\Bigl(N{\log\log q\o \log
q}\Bigr). \eqno (2.2)
$$
The contribution of terms in which $n$ is divisible by $p^2$ for
some $p$ under the condition $X<p\leq Y$ is
$$
\ll \sum_{X<p\leq Y}{N\o p^2}\ll {N\o \log^2 q}.
$$

Denote by $B_{r}$ the set of $n\le N$ not dividing by $p^{2}$ for
some $p\in (X,Y]$ and having exactly $r$ distinct prime factors in
$(X,Y]$. Then
$$
\sum_{\substack{n\leq N\\ n\in R\\ (n,\,q)=1}}f(n)\chi(n+a)
=\sum_{1\leq r\ll \log N}\sum_{\substack{n\leq N\\ n\in B_r\\
(n,\,q)=1}}f(n)\chi(n+a)\,+\,O\biggl(\frac{N}{\log^{2}{q}}\biggr).
\eqno (2.3)
$$

Among the products $pm$ ($pm\leq N$, $(p,\,m)=1$, $X<p\leq Y$, $m\in
B_{r-1})$, the number $n(\in B_r)$ appears $r$ times. Thus
$$
\sum_{\substack{n\leq N\\ n\in B_r\\ (n,\,q)=1}}f(n)\chi(n+a)
={1\o r}\sum_{\substack{X< p\leq Y\\
(p,\,q)=1}}f(p)\sum_{\substack{m\in B_{r-1}\\ mp\leq N\\ (m,\,p)=1\\
(m,\,q)=1}}f(m)\chi(pm+a).
$$
If we get rid of the condition $(m,\,p)=1$, then the corresponding
error is
$$
\ll{1\o r}\sum_{X< p\leq Y}\sum_{\substack{m\leq{N\o p}\\ p|\,m}}1
\ll{1\o r}\sum_{X< p\leq Y}{N\o p^2}\ll {1\o r}\sum_{X<p}{N\o
p^2}\ll{N\o rX}.
$$
Therefore we have
$$
\sum_{\substack{n\leq N\\ n\in B_r\\ (n,\,q)=1}}f(n)\chi(n+a)
={1\o r}\sum_{\substack{X< p\leq Y\\
(p,\,q)=1}}f(p)\sum_{\substack{m\in B_{r-1}\\ mp\leq N\\
(m,\,q)=1}}f(m)\chi(pm+a)+O\Bigl({N\o rX}\Bigr).
$$

We divide $\{X<p\leq Y\}$ into $O(\log Y)$ intervals of the form
$\{Q<p\leq Q_1\}$, where $X\leq Q\leq Q_1\leq Y$, and write
\begin{align*}
S_r(Q,\,Q_1)&=\sum_{\substack{Q< p\leq Q_1\\ (p,\,q)=1}}f(p)
\sum_{\substack{m\in B_{r-1}\\ mp\leq N\\ (m,\,q)=1}}
f(m)\chi(pm+a)\\
&=\sum_{\substack{m\in B_{r-1}\\ m\leq {N\o Q}\\ (m,\,q)=1}}f(m)
\sum_{\substack{Q< p\leq Q_1\\ p\leq {N\o m}\\
(p,\,q)=1}}f(p)\chi(pm+a).
\end{align*}
We have
$$
|S_r(Q,\,Q_1)|\leq\sum_{m\leq {N\o Q}}\Bigl|\sum_{\substack{Q< p\leq
Q_1\\ p\leq {N\o m}\\ (p,\,q)=1}}f(p)\chi(pm+a)\Bigr|.
$$
We shall estimate the sum
$$
W_1=\sum_{m\leq {N\o Q}}\Bigl|\sum_{\substack{Q< p\leq Q_1\\ p\leq
{N\o m}\\ (p,\,q)=1}}f(p)\chi(pm+a)\Bigr|. \eqno(2.4)
$$

In the following, we need

{\bf Lemma 2}. \emph{Let $q$ be a prime number, $\chi$ be a
non-principal Dirichlet character modulo $q,\,(a,\,q)=1$. Then for
two primes $p_1,\,p_2$ with $(p_1p_2,\,q)=1,\,p_1\not\equiv p_2\,
({\rm mod}\,q)$, we have}
$$
\sum_{U<m\leq V}\chi\Bigl({p_1m+a\o p_2m+a}\Bigr)\ll {V-U\o
\sqrt{q}}+\sqrt{q}\log q.
$$
Here we write ${1\o n}$ as the multiplicative inverse of $n$ such
that ${1\o n}\cdot n\equiv 1\,({\rm mod}\,q)$ and appoint ${1\o
0}=0$.

This is Lemma 4 in \cite{Gong_Jia_2015} which is produced by Lemma
18 in \cite{Bourgain_Garaev_Konyagin_Shparlinski_2012}.

By Cauchy's inequality, we have
$$
|W_1|^2\leq {N\o Q}\sum_{m\leq {N\o Q}}\Bigl|\sum_{\substack{Q< p\leq Q_1\\
p\leq {N\o m}\\ (p,\,q)=1}}f(p)\chi(pm+a)\Bigr|^2.
$$
Then
\begin{align*}
W_2&=\sum_{m\leq {N\o Q}}\Bigl|\sum_{\substack{Q< p\leq Q_1\\
p\leq {N\o m}\\ (p,\,q)=1}}f(p)\chi(pm+a)\Bigr|^2\\
&=\sum_{m\leq {N\o Q}}\sum_{\substack{Q< p_1\leq Q_1\\
p_1\leq {N\o m}\\ (p_1,\,q)=1}}\sum_{\substack{Q< p_2\leq Q_1\\
p_2\leq {N\o m}\\ (p_2,\,q)=1}}f(p_1)\overline{f(p_2)}
\chi\Bigl({p_1m+a\o p_2m+a}\Bigr)\\
&=\sum_{\substack{Q< p_1\leq Q_1\\ (p_1,\,q)=1}}\sum_{\substack{Q<
p_2\leq Q_1\\ (p_2,\,q)=1}}f(p_1)\overline{f(p_2)}
\sum_{\substack{m\leq {N\o Q}\\ m\leq{N\o
\max(p_1,\,p_2)}}}\chi\Bigl( {p_1m+a\o p_2m+a}\Bigr)\\
&\ll\sum_{\substack{Q< p_1\leq Q_1\\
(p_1,\,q)=1}}\sum_{\substack{p_1 \leq p_2\leq Q_1\\
(p_2,\,q)=1}}\Bigl|\sum_{m\leq{N\o p_2}}\chi\Bigl( {p_1m+a\o
p_2m+a}\Bigr)\Bigr|\\
&\ll\sum_{Q< p_1\leq Q_1}{N\o Q}+\sum_{\substack{Q< p_1\leq Q_1\\
(p_1,\,q)=1}}\sum_{\substack{p_1<p_2\leq Q_1\\
(p_2,\,q)=1}}\Bigl|\sum_{m\leq{N\o p_2}}\chi\Bigl({p_1m+a\o
p_2m+a}\Bigr)\Bigr|.
\end{align*}

We note that if $p_1,\,p_2\leq Y=q^{\e\o 4}$, then $p_1\ne p_2$ is
equivalent to $p_1\not\equiv p_2\,({\rm mod}\,q)$. By Lemma 2,
\begin{align*}
W_2&\ll N+\sum_{Q< p_1\leq Q_1}\sum_{Q<p_2\leq Q_1}\Bigl({N\o
\sqrt{q}Q}+\sqrt{q}\log q\Bigr)\\
&\ll N+{NQ\o \sqrt{q}}+Q^2\sqrt{q}\log q\\
&\ll N.
\end{align*}
Hence,
$$
|S_r(Q,\,Q_1)|\leq W_1\ll{N\o \sqrt{Q}}.
$$

Taking
$$
Q=2^kX,\qquad k=0,\,1,\,\ldots,\,k_0,\qquad 2^{k_0-1}X\leq Y<
2^{k_0}X,
$$
we get
$$
\sum_{\substack{X< p\leq Y\\ (p,\,q)=1}}f(p)\sum_{\substack{m \in
B_{r-1}\\ mp\leq N\\ (m,\,q)=1}}f(m)\chi(pm+a)\ll
N\sum_{k=0}^{k_0}{1\o 2^{k\o 2}\sqrt{X}}\ll {N\o \log q}.
$$
It follows that
$$
\sum_{\substack{n\leq N\\ n\in B_r\\ (n,\,q)=1}}f(n)\chi(n+a)\ll{N\o
r\log q}.
$$
Then
$$
\sum_{\substack{n\leq N\\ n\in R\\ (n,\,q)=1}}f(n)\chi(n+a)
\ll\sum_{1\leq r\ll \log N}{N\o r\log q}\ll N{\log\log q\o \log q}.
$$
Thus
$$
\sum_{n\leq N}f(n)\chi(n+a)\ll N{\log\log q\o \log q}.
$$

We have proved (1.2) for sufficiently large prime number $q$. For
bounded $q$, this is a trivial result. The proof of (1.3)
follows the same lines, but we need to replace Lemma 2 in the above
arguments by the following

{\bf Lemma 3}. \emph{Let $q$ be a prime number, $\chi$ be a
non-principal Dirichlet character modulo $q$. Assume that $t\geq
2,\,(b_1,\,q)=\cdots=(b_t,\,q)=1$, $b_1,\,\ldots,\,b_t$ are distinct
integers modulo $q$ and that primes $p_1,\,p_2$ satisfy
$(p_1p_2,\,q)=1$, $p_1\not\equiv p_2\,({\rm mod}\,q)$. Then if
$p_2\not\equiv \overline{b_i}b_j p_1\,({\rm mod}\,q)$ $(1\leq
i,\,j\leq t)$, we have
$$
\sum_{U<m\leq V}\chi\Bigl({(p_1m+b_1)\cdots(p_1m+b_t)\o (p_2m
+b_1)\cdots(p_2m+b_t)}\Bigr)\ll {V-U\o \sqrt{q}}+\sqrt{q}\log q.
$$
}

This is Lemma 6 in \cite{Gong_Jia_2015} which is produced by Lemma
17 in [1]. Lemma 17 in [1] is due to A.~Weil \cite[Appendix V,
Example 12]{Weil_1974}, which states that
$$
\biggl|\sum_{x=1}^{q}\chi_{1}(x+c_{1})\cdots \chi_{r}(x+c_{r})
e^{2\pi i\,\frac{f(x)}{q}}\biggr| \le\,(r+d)\sqrt{q}.
$$
Here $\chi_{1},\ldots,\chi_{r}$ are any Dirichlet characters modulo
$q$ such that at least one of them is non-principal, $f\in
\mathbb{F}_{q}[x]$ is any polynomial of degree $d\ge 0$ and
$c_{1},\ldots, c_{r}$ are distinct integers modulo $q$.

We should discuss in two cases respectively. In the first case,
there is no $a_i$ in (1.3) is equivalent to 0 modulo $q$. In the
second case, one of $a_i$'s is equivalent to 0 modulo $q$, which is
reduced to the case of $t-1$ by writing $f(n)\chi(n)=f_1(n)$. Then
the proof of (1.3) follows.

So far the proof of Theorem is complete.

\vskip.3in
\noindent{\bf Acknowledgements}

The first author is supported by the National Natural Science
Foundation of China (Grant No. 11671119). The second author is
supported by the National Key Basic Research Program of China
(Project No. 2013CB834202) and the National Natural Science
Foundation of China (Grant No. 11371344 and Grant No. 11321101). The
work of the last author is supported by Russian Science Foundation
(Grant No. 14-11-00433) and performed in Steklov Mathematical
Institute of Russian Academy of Sciences.


\bigskip

Ke Gong

Department of Mathematics, Henan University, Kaifeng, Henan 475004,
P. R. China

E-mail: {\tt kg@henu.edu.cn}

\

Chaohua Jia

Institute of Mathematics, Academia Sinica, Beijing 100190, P. R.
China

Hua Loo-Keng Key Laboratory of Mathematics, Chinese Academy of
Sciences, Beijing 100190, P. R. China

E-mail: {\tt jiach@math.ac.cn}

\

Maxim Aleksandrovich Korolev

Steklov Mathematical Institute of Russian Academy of Sciences

119991, Russia, Moscow, Gubkina Str., 8

E-mail: {\tt korolevma@mi.ras.ru}

\end{document}